\documentstyle[12pt]{article}
\def~{\unskip\nobreak\ }
\newcount\hh
\newcount\mm
\mm=\time
\hh=\time
\divide\hh by 60
\divide\mm by 60
\multiply\mm by 60
\mm=-\mm
\advance\mm by \time
\def\hhmm{\number\hh:\ifnum\mm<10{}0\fi\number\mm} %
\makeatother
\title{Random Banach spaces.\\
The Limitations of the Method}
\author{\small Piotr Mankiewicz 
\thanks{Part of this research was done while this author has been visiting 
the Case Western Reserve University. Supported in part by a grant from KBN (Poland).}\\
\small Institute of Mathematics\\
\small Polish Academy of Sciences\\
\small Warsaw
\and
\small Stanis{\l}aw J. Szarek
\thanks{Supported in part by a grant from the National Science Foundation (U.S.A.).}\\
\small Case Western Reserve University\\
\small Cleveland}
\date{}
\newtheorem{fact}{Fact}
\newtheorem{thm}{Theorem}[section]
\newtheorem{prop}[thm]{Proposition}
\newtheorem{lemma}[thm]{Lemma}
\newtheorem{cor}[thm]{Corollary}

\newcommand{\pend}{\bigskip\hfill\(\Box\)} 
\newcommand{\Rn}[1]{\mbox{{\it I\kern -0.25emR}$\sp {\,{#1}}$}}
\newcommand{\Nn}[1]{\mbox{{\it I\kern -0.25emN}$\sp {\,{#1}}$}}
\newcommand{\En}{\mbox{{\it I\kern -0.25emE}}} 

\begin{document}

\maketitle

We shall study the properties of typical $n$--dimensional subspaces of
  $l^N_{\infty} =
(\Rn{N}, \| \cdot \|_{\infty})$, or equivalently, typical $n$--dimensional
quotients of $l^N_1 = (\Rn{N}, \| \cdot \|_1)$, where the meaning what is
 typical and what 
is not is defined in terms of the Haar measure $\mu_{n,N}$ on the Grassmann 
manifold $G_{n,N}$ of 
all $n$--dimensional subspaces of $\Rn{N}$. \par

In [Gl.2], Gluskin proved that a ``typical'' $n$--dimensional subspace $E$
 of $l^{n^2}_{\infty}$ 
enjoys the property $$
\|P\| \geq \frac{ck}{\sqrt{n\log n}},
$$
for every projection $P:E \rightarrow E$, with $\min\{{\rm rank}\, P,
{\rm rank}\, (Id-P)\} = k$, 
where $c$ is a numerical constant. In particular, if $k\geq n^{\alpha}$, 
$\alpha>\frac{1}{2}$, 
then no projection P on E with both ${\rm rank}\, P$ and ${\rm corank}\, 
P$ greater than $k$ can 
be ``well'' bounded.
Several other results,[Sz.1],[Sz.2],[Ma.1],[Ma.2] showed that a ``typical''
proportional (i.e. $\dim E \approx \beta N$ for some ``fixed'' $\beta \in 
(0,1))$ 
subspace $E$ of $l^N_{\infty}$ has the property that every ``well'' bounded 
operator on $E$ is 
indeed a
``small'' perturbation of a multiple of the identity $\lambda {\rm Id}_E$. 
However, the estimates on the 
distance between $T$ and $\lambda {\rm Id}_E$ have been done in terms of the 
geometry of $\Rn{N}$ rather 
than $E$ itself. In this note, we obtain the estimates on the distance 
between $T$ and $\lambda 
{\rm Id}_E$ in intrinsic terms of the geometry of $E$, namely, in terms of 
the Gelfand numbers of 
$T - \lambda {\rm Id}_E$ (Sections 2 and 3). On the other hand, we show in 
Section 4, that if 
$k \leq n^{1/2}$ then a ``typical'' $n$--dimensional subspace $E$ of 
$l^N_{\infty}$ (for 
any $N\geq n$) contains a $k$--dimensional well-complemented subspace $G$ 
isomorphic to $l^k_p$ 
with either $p=2$ or $p=\infty$ and therefore admits operators which are 
``fairly'' far away from the line $\{\lambda {\rm Id}_E\}_{\lambda \in \Rn{}}$.
\par

We shall employ the standard notation of local theory of Banach spaces as 
used in  e.g. [F-L-M].
For basics on multivariate Gaussian random variables the reader is 
referreded  to [T].

\section{{\protect\large\bf Generic subspaces. Equivalent Gaussian approach}}
Let $P_n$ for $n \in \Nn{}$ be a sequence of properties of $n$--dimensional 
Banach spaces and let $f: \Nn{} \rightarrow \Nn{}$ be an increasing function. 
We shall say that $P=\{P_n\}_{n \in \Nn{}}$ is a generic property of 
$n$--dimensional subspaces of $l^N_{\infty}$ where $N= f(n)$ iff for 
every $n \in \Nn{}$ we have $$
\mu_{n,N} \{E \in G_{n,N} \,| \, E \,\, {\rm satisfies} \,\, P_n \} 
\geq 1 - \varepsilon^n ,
$$
for some $\varepsilon \in (0,1)$, $\varepsilon$ independent of $n$. 
In the sequel, we shall say that a generic $n$--dimensional subspace of 
$l^N_{\infty}$ has a property $P=\{P_n\}_{n \in \Nn{}}$ rather then that 
the property $P$ is a generic one. Since the dual of an $n$--dimensional 
subspace $E$ of $l^N_{\infty}$ is the quotient $F=l^N_1 /E^{\perp}$, where
$l^N_1 = (\Rn{N}, \|\,\cdot\,\|)$
and $E^{\perp}$ is the orthogonal complement of $E$ in $\Rn{N}$, the notions 
of generic properties of quotients of $l^N_1$ and generic subspaces of 
$l^N_1$ can be defined. E.g. a generic quotient of $l^N_1$ satisfies a 
property $P$ iff the corresponding (via duality) geometric generic subspace 
of $l^N_{\infty}$ satisfies the dual property $P^*$.\par 

In the context of quotients it is more convenient to consider an equivalent 
approach.\par

Let $g$ be an $\Rn{n}$--valued Gaussian vector distributed according to the 
$N(0, n^{-1}{\rm Id}_{\Rn{n}})$ law (i.e., the covariance matrix of $g$ 
is $n^{-1}{\rm Id}_{\Rn{n}}$, or the coordinates of $g$ are i.i.d. 
$N(0,1/n)$ random variables). Let $g_1,g_2,\ldots,g_N$ be independent 
copies of $g$ and $\Gamma$-- an $n \times N$ matrix whose columns are 
$g_1,g_2,\ldots,g_N$; alternatively, $\Gamma=\Gamma(\omega)$ can be 
described as a Gaussian $n \times N$ matrix with i.i.d. $N(0,1/n)$ entries. 
If we think of $\Gamma$ as of a (random) linear map from $\Rn{N}$ to $\Rn{n}$, 
$l_1^N/{\ker \Gamma}$ is a random quotient of $l^N_1$. By the rotational 
invariance of the Gaussian measure, this model is ``measure theoretically'' 
equivalent to the one described at the begining of this section and based 
on the Haar measure on the Grassmannian. Still equivalently, we may 
consider the random norm on $\Rn{n}$, whose unit ball is a random 
absolute convex body
$$
B=B(\omega)= {\rm absconv}\,\{g_j\,:\, j=1,2,\ldots,N\} =\Gamma(B^N_1)\,\, ,
$$
where $B^N_1$ is the unit ball of $l^N_1$.\par 

In the sequel we will need some basic facts about Gaussian vectors, 
Gaussian matrices and ``Gaussian bodies''.The first lemma is an elementary 
consequence of the formulae for Gaussian density in $\Rn{d}$.\par 

\begin{lemma}
\label{A}
If $g$ is a Gaussian random variable with distribution\\ \medskip
$N(0,d^{-1}{\rm Id}_{\Rn{d}})$, then\\
\medskip
\mbox{} \hspace{1cm} {\rm (i)} $\En\,\|\,g\,\|^2_2=1$\\ \medskip
\mbox{} \hspace{1cm} {\rm (ii)}
${\bf P}\left(\|\,g\,\|_2 \geq \lambda \right) \leq \exp{(-d\lambda^2/8)} 
\mbox{} \hspace{.7cm} {\rm for}\hspace{.3cm} \lambda \geq 2$\\ \medskip
\mbox{} \hspace{1cm} {\rm (iii)}
${\bf P} \left(\|\,g\,\|_2 \leq t \right) \leq 
\left(t{\rm e}^{1/2}\right)^d$\\
\medskip
\mbox{} \hspace{1cm} {\rm (iv)}
${\bf P} \left(\frac{1}{2} \leq \|\,g\,\|_2 \leq 2 \right) \geq 
1 - \exp{(-c_0 d)}$\,,\\
where $c_0$ is a universal constant.
\end{lemma}

The next lemma can be derived from the first one using a standard 
$\epsilon$--net argument (cf. 
[Sz.3], Lemma 2.8). Another, more precise argument (giving e.g. 
$c>\frac{1}{4}$ and $C<2$) can be 
found in [Si].\par

\begin{lemma}
\label{B}
Let $ k\leq \frac{1}{2}N$ and $\Lambda$ be an $N \times k$ matrix with all 
i.i.d. Gaussian $N(0,1)$ entries. Then $${\bf P} (\, c\|\,x\,\|_2 \leq 
k^{-1/2}\|\,\Lambda x\,\|_2 \leq C\|\,x\,\|_2 \,\, \mbox{{\rm for every}} 
\, x \in \Rn{k}) \geq 1 - \exp{(-c_1 d)}\, ,
$$
where $c,C,c_1$
are universal constants.
\end{lemma}

We have an immediate\\

\begin{cor}
\label{C}
If $k$ and $N$ are as in Lemma {\ref B}, $g_1,g_2,\ldots, g_N$ are i.i.d. 
Gaussian ran\-dom var\-ia\-bles with dis\-tri\-bu\-tion $N(0,k^{-1}
{\rm Id}_{\Rn{k}})$ and $B=B(\omega)= {\rm absconv}\, 
\{g_1,g_2,\ldots ,g_N\}$, then
$$
{\bf P} (B \supset ck^{-1/2}D) \geq 1 - \exp{(-c_1 d)}\, , $$
where $D$ stands for the Euclidean unit ball in $\Rn{k}$. \end{cor}

The last lemma gives more precise information about the random bodies 
$B(\omega) \subset \Rn{k}$ (cf. [Gl.3],[Gl.4]). 

\begin{lemma}
\label{D}
If $g_1,g_2,\ldots,g_N$, $B$ are as in Corollary {\ref C} and \medskip
$2^k\geq N \geq 2k$, then\\
\medskip
\mbox{} \hspace{1cm} {\rm (i)} ${\bf P}
\left( B \supset c'\sqrt{\frac{\log{N/k}}{k}}D \right) \geq 
1 - \exp{(-c_2k)}$\\
\medskip
\mbox{} \hspace{1cm} {\rm (ii)} ${\bf P} 
\left( ({\rm vol}\,B\,/ {\rm vol}\,D\,)^{1/k} 
\leq C'\sqrt{\frac{\log{N/k}}{k}}\right) 
\geq 1 - \exp{(-c_2k)}\, ,$\\ where $c,C', c_2 >0$ are universal 
constants. \end{lemma}

\begin{rem}
The volume estimate from Lemma {\ref D} (ii) actually holds for any 
$B= {\rm absconv}\,
\{x_1,x_2,\ldots,x_N\}$ as long as we control 
$$ \max{\{\|\,x_j\,\|_2\,:\, j= 1,2,\ldots,N \}},$$
which in our case we do by Lemma {\ref A} (iv) 
(see [Ca-P] or [Gl.4]). \end{rem}

\section{{\protect\large\bf The proportional case.}} 

Our starting point is the following result ([Ma.2], Proposition 2.3). 

\begin{thm}
\label{baza}
There is a numerical constant $c>0$ such that for every $n \geq 2$ there 
is a norm $\| \cdot 
\| \sb {X \sb n }$ on $\Rn{n}$ such that\\
\mbox{} \hspace{1cm} {\rm (i)} $ X \sb n = (\Rn{n}, \| \cdot \|)$ is 
isometrically isomorphic to a quotient of $l \sp {2n} \sb 1$,\\ 
\mbox{} \hspace{1cm} {\rm (ii)} $ \|\,x\,\| \sb 2 \leq \|\,x\,\| 
\sb {X \sb n} \leq \|\,x\,\| \sb 1$ \,\, \,\, for every $x \in \Rn{n}$,
\\ \mbox{} \hspace{1cm} {\rm (iii)} for every $ T \in L(\Rn{n})$ 
there are $\lambda \sb T \in \Rn{}$ , $ V \sb T \in L(\Rn{n})$ and 
a linear subspace $E_T \subset \Rn{n} $ with ${\rm dim}\,E \sb T > 
\frac{7n}{8}$, such that\\
\mbox{} \hspace{1.5cm} {\rm a)} $ V \sb T = T + \lambda {\rm Id} 
\sb {\Rn{n}}$ ,\\
\mbox{} \hspace{1.5cm} {\rm b)} $ |\,\lambda \sb T \,| 
\leq c\|\,T\,\| \sb {X \sb n}$ ,\\
\mbox{} \hspace{1.5cm} {\rm c)} $ \|\, V\sb T\,| E \sb T\,\| 
\leq cn \sp {-\frac{1}{2}} \|\,T\,\| \sb {X \sb n}$ . \end{thm}

In fact, Theorem 2.1 holds for ``sort of'' generic $n$-dimensional 
quotients of $l \sp {2n} \sb 1$
({\em cf.} [Ma.1]).  However in order to adapt the above result to 
our present setting, we
need to make a couple of observations.  First, the condition (ii), 
which is not crucial for our
purposes, 
has to be superseded by the
properties listed in Corollary 1.3 and Lemma 1.4 (the condition of 
type (ii) may be, moreover,
achieved, up to a universal constant and restricted to the span of, 
say, the first
$n/2$ unit vectors $e_1,e_2, \ldots, e_{n/2}$; see [Ma.-T]). 
Next (and more importantly), we have to point out
that in all constructions leading to Theorem 2.1--like statements 
([G1.1], [G1.2], [Ma.1], [Ma.2],
[Ma.-T], [Sz.1], [Sz.2] etc.), ``generic'' had a somewhat different 
(and slightly less
natural) meaning.  We take here an opportunity to present a remark 
which rectifies this problem.
What happens is that when ensuring the condition (iii) from Theorem 
2.1, we need to work with
the set

$${\cal T} = \{ T \in L( \Rn{n}) : \| T \|_{X_n} \leq 1 \}$$
or, more specifically, with nets of ${\cal T}$ in the $l \sp {N} 
\sb 1$ metric.  Now we
have $(hs( \cdot )$ is the Hilbert-Schmidt norm)

\medskip
{\em If $X,Y$ are quotients of $l \sp {N} \sb 1$, endowed with the 
canonical inner product and
$T : X \rightarrow Y$ verifies $\| T : X \rightarrow Y \| \leq 1$, 
then $hs(T) \leq
N^{1/2}$.}
\medskip

The above statement is shown by estimating 
$hs(T) = hs(T^*) = \pi_2 (T^*)$ by the $\pi_2$-norm of
(a restriction of) the formal identity ${\rm Id} : l_{\infty}^N 
\rightarrow \l_2^N$, which is
$N^{1/2}$.  If $N = 2n$, it can be shown in a standard way that 
${\cal T}' = \{ T \in
L( \Rn{n}) : hs (T) \leq N^{1/2} \}$ admits a $\delta$-net in the 
$\| \cdot \|_{l_2^n}$
metric which is of cardinality not exceeding $(C/ \delta)^{n^2}$, 
where $C$ is a numerical constant, 
and this can be easily incorporated into existing proofs of 
Theorem 2.1-like statements.
Unfortunately, we do not see how to handle in the same ``unified'' 
way e.g. the casse considered 
in Theorem 3.1.

We now prove the following.

\begin{thm}
\label{n/2}
There is a numerical constant $K > 0$ such that if $X \sb n$ is the 
quotient space from Theorem \ref{baza} then, for every operator $T 
\in L(\Rn{n})$, we have

$$
\inf \{ c \sb {\frac{n}{2}} (T - \lambda{\rm Id} 
\sb {\Rn{n}}) \, |\, \lambda \in \Rn{} \, \}
\leq K n^{- \frac{1}{2}} \| T \|_{X_n}. $$
\end{thm}
We recall here that, for an operator $u: X \rightarrow Y$, the $k$th 
Gelfand number of $u$ is
defined by
$$c_k(u) = \inf \{ \| u_{|Z} \| : Z \subset X, \, \mbox{ codim } Z < k \}.$$
Because of a well-known duality relation between the Gelfand numbers and 
the so-called 
``Kolmogorov numbers'' $d_k ( \cdot )$  (namely $d_k(u^*) = c_k (u))$, 
one can also state
our results in terms of the latter ones.  Note a slight abuse of notation; 
in the above and in what follows we pretend that $\frac{n}{2}$ 
and similar expressions are integers.

Observe that Theorem 2.2 is the best possible.  This follows either 
from Corollary 4.3 below or from
the fact that, for an $n$-dimensional ``generic'' quotient of $l_1^N$, 
(for any $N > n)$,
and for a ``generic'' element of $O(n)$,  the left hand side of the 
inequality in Theorem 2.2
is of order $n$ while the right hand side--at most of order $n^{1/2}$.  
In fact, any 
$n$-dimensional normed space can be represented on $\Rn{n}                                                   $ so that the last remark holds. 

{\bf Proof of Theorem 2.2}
Obviously, it is enough to prove the theorem for every operator $T 
\in L(\Rn{n})$ satisfying
\begin{equation}
\label{-1}
\inf \{ c \sb {\frac{n}{2}}
(T - \lambda{\rm Id} \sb {\Rn{n}}) \, | \lambda \in \Rn{} \, \} = 1 .
\end{equation}
To this end, fix such an operator T. It is well known, [Gl.1], [Sz.2], that
$$
{\rm vol}\, (B \sb {X \sb n} ) \leq \left ( \frac{c \sb 1}{n} \right) \sp n $$
and
$$
B \sb {X \sb n} \supset \frac{1}{\sqrt n} B \sp 2 \sb n . $$
Hence, by [Sz.-T], we infer that there exists an 
$\frac{3n}{4}$\,-\,dimensional subspace, say $E$, of $\Rn{n}$ such that 

\begin{equation}
\label{1}
B \sb {X \sb n} \cap E \, \subset \, \frac{c \sb 2}{\sqrt n} B \sb E 
\sp 2. \end{equation}

{\bf Claim}. For every $\lambda \in \Rn{n}$ the operator \begin{equation}
\label{2}
T \sb {\lambda}\,|E = (T - \lambda {\rm Id} \sb {\Rn{n}})\,|E \end{equation}
has at least $n/4$\, $s$\,-\,numbers greater than or equal 
to $1/c \sb 2$ .\par Indeed, if this was not the case then, 
by (\ref{1}), we would have for some $\lambda \sb 0 \in \Rn{}$
\begin{eqnarray}
\frac{1}{c \sb 2} > \|\, T \sb {\lambda \sb 0}\,|E \sb 0 : (E \sb 0, n 
\sp {-\frac{1}{2}} \|\,.\,\| \sb 2\,) \rightarrow \,(\Rn{n}, n 
\sp {-\frac{1}{2}} \|\,.\,\| \sb 2)\,\| \nonumber
\geq \\
\frac{1}{c \sb 2} \|\, T \sb {\lambda \sb 0}\,|E \sb 0 : 
(E \sb 0, \|\,.\,\| \sb {X \sb n}\,) \rightarrow \,(\Rn{n}, \|\,.\,\| 
\sb {X \sb n})\,\| , \label{3}
\end{eqnarray}
where $E \sb 0$ is an $\frac{n}{2}$\,-\,dimensional subspace such that 
$\|\,T \sb {\lambda \sb 0}\,|E \sb 0 \,\| \sb 2 < 1/c \sb2$. 
But (\ref{3}) implies that
$$
c \sb {\frac{n}{2}}(T - \lambda \sb 0 {\rm Id} \sb {\Rn{n}}) < 1 , $$
a contradiction with (\ref{-1}), which concludes the proof of the claim. 

In particular, the Claim yields that for every $\lambda \in \Rn{}$ 
the operator $T - \lambda {\rm Id} \sb {\Rn{n}}$ has at least 
$n/4$ $s$\,-\,numbers greater then or equal to $1/c \sb 2$, 
which means that for every subspace $F \subset \Rn{n}$ with $\dim F 
\geq \frac{7n}{8}$ we have
$$
\|\,T - \lambda {\rm Id} \sb {\Rn{n}}\,|F\,\| \sb 2 \geq \frac{1}{c \sb 2} .
$$

Hence, by Theorem \ref{baza} (iii), c), we infer that $$ 1 
\leq c \sb 2 cn \sp {- \frac{1}{2}} \|\,T\,\| \sb {X \sb n} $$
which implies $\|\,T\,\| \sb {X \sb n} \geq Kn \sp {\frac{1}{2}}$, 
where $K = (cc \sb 2) \sp {-1}$. This concludes the proof of the theorem.
\pend

\section{{\protect\large\bf The $l\sp 1$ - type estimates.}} 

For a Banach space $X=(\Rn{n},\|\,.\,\|_B)$ we set $$
M^*_B = \int_{S^{n-1}} \|\,x\,\|^*_B d\mu(x), $$
where $d\mu$ stands for the normalized Lebesque measure on the unit 
sphere $S^{n-1}$ and $\|\,.\,\|^*_B$ 
denotes the dual norm to $\|\,.\,\|_B$; this is the ``average width'' 
of $B$.\par

In the sequell we shall need the following fact which can be found in 
[P--T], Theorem 1 ({\em cf.}\, [Pi.2 ], Theorem 1.3). 

\begin{fact}
\label{Pisier}
There exist a numeric constant
$C > 1$ such that for every symmetric convex body $B \subset \Rn{n}$ 
and for every $k = 1,2,\ldots,n-1$
there exists a subspace $E \sb k \subset \Rn{n}$ with ${\rm codim}\, 
E \sb k = k$ such that
$$
B\cap E \sb k \subset C M \sp * \sb B
\sqrt {\frac{n}{k}} D \cap E \sb k.
$$
\end{fact}

Recall that for an operator $T \in L(\Rn{n})$ we say that 
$T \in M \sb n (\alpha,\beta)$, where $\alpha,\beta > 0$ iff 
there is a linear subspace $F \subset \Rn{n}$ with 
${\rm dim}\, F \geq \alpha$ such that $$
\|P \sb {F \sp {\perp}}Tx\| \sb 2 \geq 
\beta \|x\| \sb 2 \,\, \mbox{for every } x \in F ,
$$
(where $P \sb {F \sp {\perp}}$ denotes the orthogonal projection onto 
$F \sp {\perp}$). Also, for $\gamma > 0$, we denote $$
\tilde{M} \sb n (\gamma) = \bigcup\sb{k=1}\sp{n/2} M \sb n (k,\gamma /k). $$

The following fact has been proved in [Sz.2] 

\begin{fact}
\label{mix}
A generic n\,-\,dimensional quotient $X_n$ of $l \sp {n \sp 2} \sb 1$ 
enjoys the property $$
\|T\| \sb {X_n} \geq \frac{c \sb 1 \gamma}{\sqrt {n \log n}} $$
for every $T \in \tilde
{M} \sb n (\gamma)$.
\end{fact}

\medskip\noindent From Fact I and Fact II we deduce

\begin{thm}
\label{thm.2}
A generic n\,-\,dimensional quotient $X_n$ of $l \sp {n \sp 2} \sb 1$ 
has the property that for every $T \in L(\Rn{n})$ we have
$$
\inf {\{\sum \sb {i=1} \sp n c \sb i 
(T - \lambda {\rm Id} \sb {\Rn{n}}) \,\,|\,\,\lambda \in \Rn{}\}} 
\leq c n \sp {2/3} \sqrt {\log \sp 3 n}\,\,\|T\| \sb {X \sb n} ,
$$
where $c$ is a numerical constant.
\end{thm}

\noindent {\bf Remark.}  It is imaginable that one could strengthen 
the above inequality to get
$O(n^{1/2}) \|T \|_{X_n}$ on the right hand side; 
$\circ (n^{1/2}) \| T \|_{X_n}$ is
impossible, see the comments following Theorem 2.2.

\medskip\noindent {\bf Proof}
In order to simplify the notation we shall assume that 
$n = 2 \sp k$ for some $k \in {\bf N}$.
Let $X_n$ be a generic n--dimensional quotient of $l \sp 1 \sb {n \sp 2}$.
Obviously, it suffices
to prove that
\begin{equation}
\label{!1}
\sum \sb {i=1} \sp n c \sb i (T)
\leq c n \sp {2/3} \sqrt{\log \sp 3 n}\,\,\|T\| \sb {X_n},
\end{equation}
for every $T \in L(\Rn{n})$ satisfying ${\rm tr}\,T = 0$. 
To this end, fix $T \in L(\Rn{n})$ such that ${\rm tr}\,T = 0$ and
\begin{equation}
\label{!2}
\sum \sb {i=1} \sp n c \sb i (T) = n .
\end{equation}
Then, there is $i \leq k$ such that
\begin{equation}
\label{!3}
2 \sp {i-1} c \sb {2 \sp i} \geq \frac{n}{\log n}. \end{equation}
If $2 \sp i \leq n \sp {2/3}$ then we have $$
\|T\|
\sb {X_n} \geq c \sb {2 \sp i} \geq \frac{2n \sp {1/3}}{\log n} \,, $$
which combined with (\ref{!2}), by a standard homogenuity argument, 
yields (\ref{!1}) and we are done.

Thus, assume that $2 \sp i > n \sp {2/3}$. It is well known ({\em cf.} 
Corollary \ref {C}) that $n \sp {- 1/2} D \subset B \sb {X(E)}$
and that $M \sp * \sb {X(E)
} \leq
c \sb 2 \sqrt{n \sp {-1} \log n}$, ({\em see e.g.} [F-L-M]). 
Thus, applying Fact \ref{Pisier} we infer that there exists a linear 
subspace $F \sb {2 \sp {i-1}} \subset \Rn{n}$ 
with ${\rm codim}\,F \sb {2 \sp {i-1}} = 2 \sp {i-1}$ such that 
\begin{equation}
\label{!4}
\frac{1}{\sqrt n} B \sp 2 \sb n \cap F \sb {2 \sp {i-1}} 
\subset B \sb {X \sb n} \cap F \sb {2 \sp {i-1}} \subset 
c \sb 2 C \sqrt {\frac {\log n}{2 \sp {i-1}}} B \sp 2 \sb n \cap F \sb 
{2 \sp {i-1}}.
\end{equation}

{\bf Claim}.
The operator $T\,|F \sb {2 \sp {i-1}}$ has at least 
$2 \sp {i-1}$ \,s\,-\,numbers greater than or equal to $$
(c \sb 2 C) \sp {-1} \sqrt {\frac{n}{2 \sp {i-1} \log \sp 3 n}}\,\, . $$
\par

Indeed, assume to the contrary that there exists a linear subspace 
$F \subset F \sb {2 \sp {i-1}}$ with ${\rm codim}\, F = 2 \sp i$ such that
\begin{equation}
\label{claim}
\|\,\,T\,|F\,:\,(F,\| \cdot \| \sb 2)
\,\rightarrow\,(\Rn{n},\,\| \cdot \| \sb 2\,)\,\| <
(c \sb 2 C) \sp {-1} \sqrt {\frac{n}{2 \sp {i-1} \log \sp 3 n}}\,\, . 
\end{equation}
Then, by (\ref{!4}) and (\ref{!3}) we have \begin{eqnarray}
(c \sb 2 C) \sp {-1} \sqrt {\frac{n}{2 \sp {i-1} \log \sp 3 n}} >
\|\,\,T\,|F\,:\,(F,\,n \sp {-1/2}\| \cdot \| \sb 2)\,\rightarrow \,
(\Rn{n},\,n \sp {-1/2}\| \cdot \| \sb 2)\,\,\| =
\nonumber \\
(c \sb 2 C) \sp {-1} \sqrt {\frac{2 \sp {i-1}}{n \log n}} \|\,\,T\,|F\,:
\,(F,\,c \sb 2 C 
\sqrt {\frac{\log n}{2 \sp {i-1}}} \| \cdot \| \sb 2)\,\rightarrow
\,(\Rn{n},\,n \sp {-1/2}\| \cdot \| \sb 2)\,\,\| \geq \nonumber \\
(c \sb 2 C) \sp {-1} \sqrt {\frac{2 \sp {i-1}}
{n \log n}} \|\,\,T\,|F\,:\,(F,\,\| \cdot 
\| \sb {X \sb n})\,\rightarrow \,
(\Rn{n},\,\| \cdot \| \sb {X \sb n})\,\,\|\geq \,\,\\ (c \sb 2 C) \sp {-1} 
\sqrt {\frac{2 \sp {i-1}}{n \log n}} c \sb {2 \sp i}(T) \geq
(c \sb 2 C) \sp {-1} \sqrt {\frac{2 \sp {i-1}}
{n \log n}} \,\,\frac{n}{2 \sp {i-1} \log n} =
\nonumber \\
(c \sb 2 C) \sp {-1} \sqrt {\frac{n}
{2 \sp {i-1} \log \sp 3 n}}\,\, , \nonumber
\end{eqnarray}
a contradiction which concludes the proof of the claim.\par 

Now, observe that if the median s\,-\,number of $T$ $$
s \sb {n/2} (T) < \frac{1}{2}(c \sb 2 C) \sp {-1} 
\sqrt {\frac{n}{2 \sp {i-1} \log n}}
$$
then, by [Ma.1], Lemma 2.6, we obtain that $$
T \in \tilde{M} \sb n \left( \frac{1}{32} 
(c \sb 2 C) \sp {-1} \sqrt {\frac{2 \sp i n}{\log n}}\, \right)\,\, , $$
while if
$$
s \sb {n/2} (T) \geq \frac{1}{2}(c \sb 2 C) \sp {-1} 
\sqrt {\frac{n}{2 \sp {i-1} \log n}}
$$
then, by [Ma.1], Theorem 3.1,
we infer that
$$
T \in \tilde{M} \sb n \left(\frac {c \sb 3}{2} 
(c \sb 2 C) \sp {-1} \sqrt {\frac{n \sp 3}{2 \sp {i-1}}}\,\right)\,\, , $$
where $c \sb 3 < 1$ is the constant from Theorem 3.1 in [Ma.1]. 
In the first case Fact II yields
$$
\|T\| \sb {X_n} \geq
\frac{c \sb 1}{32} (c \sb 2 C) \sp {-1}
\sqrt {\frac{2 \sp i}{\log n}} \geq
\frac{c \sb 1}{32} (c \sb 2 C) \sp {-1}
\sqrt{\frac{n \sp {2/3}}{\log n}}\,\, ,
$$
while in the second case, by Fact II, we get $$
\|T\| \sb {X_n} \geq
\frac{c \sb 1 c \sb 3}{2} (c \sb 2 C) \sp {-1} 
\sqrt{\frac{n}{\log \sp 3 n}}\,\, .
$$
thus, by (\ref{!2}), in both cases we obtain $$
\sum \sb {i=1} \sp n c \sb i (T) \leq
32 c \sb 2 C (c \sb 1 c \sb 3 ) \sp {-1} n \sp {2/3} 
\sqrt {\log \sp 3 n} \,\, \|T\| \sb {X_n} \, ,
$$
which proves (\ref{!1}) with $c = 32 c \sb 2 C (c \sb 1 c \sb 3 ) 
\sp {-1}$ and completes the proof of the theorem.\par \pend

\section{{\protect\large\bf The positive statements.}} 
In this section we prove several ``positive'' statements about 
existence of ``nontrivial''
operators on generic Banach spaces, which will show that the results 
of the preceeding section are ``essentially'' optimal 
(cf. Cor. \ref{elpe}). 
Results similar to some of the presented below were obtained independently 
by Gluskin [Gl.5].
The first of
these statements will also show that, for generic Banach spaces, 
the following conjecture due to 
Pisier [Pi.1] (and usually referred to as the ``dichotomy conjecture'') 
holds.\par
\vspace{.5cm}

{\em There exist $C>1$ and a function $f : \Rn{+} \rightarrow 
\Rn{+}$ with $\lim_{\lambda \to \infty} f(\lambda)= \infty$, 
such that if $E \subset l^N_{\infty}$, then there exists 
$F\subset E$, $\dim F = k \geq f(\dim E/\log N)$ verifying 
${\rm{d}}(l^k_{\infty},F) \leq C$.}\par \vspace{.5cm}

At the time of this writing the conjecture is open for general $E$. 
However, we have
\begin{prop}
\label{generic}
There exist universal constants $C,c>0$ such that if $E$ is a generic 
quotient of $l^N_1$, $\dim E = d$ (resp. a generic subspace of 
$l^N\sb{\infty}$), then there is a subspace $F \subset E$, $\dim F = k 
\geq c \min \{d\sp{\frac{1}{2}},d/\log N\}$ verifying\par 

\mbox{} \hspace{1cm} {\rm (i)} ${\rm{d}}(l^k_1,F) \leq C$ 
(resp. ${\rm{d}}(l^k\sb{\infty},F) \leq C$)\par 

\mbox{} \hspace{1cm} {\rm (ii)} $F$ is $C$--com\-ple\-men\-ted in $E$. 
\end{prop}
\medskip

\begin{prop}
\label{square}
If, in the notation of Proposition \ref{generic}, $N \geq d^2$, 
then there is $G\subset E$, $\dim G = h \geq \min\{c\log N,d\}$, 
satisfying\par

\mbox{} \hspace{1cm} {\rm (i)'} ${\rm{d}}(l^h_2,G) \leq C$\par 

\mbox{} \hspace{1cm} {\rm (ii)'} $G$ is $C$--complemented in $E$. 
\end{prop}
\medskip

\noindent
\begin{rem} \label{alpha}
If, in Proposition \ref{square}, one assumes that $N>d\sp{1+\alpha}$ 
for $\alpha \in (0,1)$, one gets $G$ which is $C/\sqrt{\alpha}$ 
complemented in $E$; for arbitrary $N$ and $d$ we get 
$C\sqrt{\frac{\log N}{\log {N/d}}}$ -- complementation. \end{rem}
\medskip

\begin{cor}
\label{elpe}
A ``generic'' $d$--dimensional quotient of $l^N_1$ contains a 
$C$--com\-ple\-men\-ted sub\-spa\-ce $C$--iso\-mor\-phic to $l^k_p$, 
$k\geq c\sqrt{d}$, either for $p=1$ or $p=2$ (resp. a ``generic'' 
subspace of $l^N\sb{\infty}$ contains such a subspace with either 
$p=\infty$ or $p=2$). \end{cor}

\noindent
{\bf Proof of Corollary \ref{elpe}} \,\, If $\log{N} < d^{\frac{1}{2}}$, 
we use Prop. 
\ref{generic} to get a $C$--complemented subsapce $C$--isomorphic to 
$l^k_1$ 
(resp. $l^k_{\infty})$. If $\log N > d^{\frac{1}{2}}$ (hence $N>d^2$), 
use Prop \ref{square} 
to get a $C$-complemented Hilbertian subspace.\par \pend

\medskip\noindent
{\bf Proof of Proposition \ref{generic}} \,\, We will prove the statement 
for a random quotient $Q$; the ``subspace'' variant follows by duality.
\par Recall that $g_j= Qe_j$, $j=1,2,\ldots,N$ are independent 
Gaussian vectors with distribution $N(0,d^{-1}Id_{\Rn{d}})$
and that the unit ball B of our random quotient is 
${\rm absconv}\{g_1,g_2,\ldots,g_N\}$. Clearly, we can assume that 
$N \leq \exp\{cd\}$. It then follows from Lemma {\ref A} (iv) that 
\begin{eqnarray}
\label{2/2}
{\bf P}(\frac{1}{2} \leq \|\,g_j\,\| \leq 2 \,\,
{\rm for} \, \, j=1,2,\ldots,N) \geq
1 - \exp (-c_1d)
\end{eqnarray}
provided that c is chosen to satisfy $c\leq c_0/2$, where $c_0$ is 
the constant from Lemma {\ref A}. Moreover, if $k \leq d/2$ and 
$ A\subset \{1,2,\ldots,N\}$ with $|A|=k$, then (cf. Lemma~{\ref B})
\begin{eqnarray}
\label{el2}
\| \, \sum_{j \in A} t_j g_j \, \| \geq c_2 
\left( \sum_{j \in A} |t_j|^2 \right)^{\frac{1}{2}} \hspace{.7cm} 
\mbox{{\rm for all choices of scalars\,}}\, (t_j)_{j \in A} 
\end{eqnarray}
with the similar probability as in (\ref{2/2}). In fact, since 
${{N}\choose{k}} <
(\frac{Ne}{k})^k$, (\ref{el2}) happens
for all such $A$ with the same estimate on the probability as in 
(\ref{2/2}) provided
$k \log{\frac{N}{k}} \leq c_3 d$. In particular, this happens if
$k \leq c \frac{d}{\log N}$ (we do not use this fact here).\par

In the next step we shall show that, for fixed 
$A \subset \{1,2,\ldots,N\}$ with $ |A| \leq k =
c\min {\{d^{\frac{1}{2}}, \frac{d}{\log N}\}}$, and 
$E= [g_j \,|\,\, j \in A ]$, we have \begin{eqnarray}
\label{fin}
{\bf P} \left(\|\, P_E g_j \,\| \leq
k^{\frac{1}{2}} \,\, | \,\, j \notin A \right) 
\geq 1 - \exp \left( -c_4\frac{d}{k}\right) \end{eqnarray}

Observe that (\ref{el2}) and (\ref{fin}) imply the conclusion of the 
Proposition \ref{generic} with $C=c_2^{-1}$. 
To this end, note that the operator $u:l^k_1 \rightarrow E$ sending 
$\{e_1,e_2,\ldots,e_k\}$ into $\{g_j \,|\,j \in A \}$ is of norm $1$, 
while $\|\,u^{-1}P_E\,\| \leq c_2^{-1}$
(notice that (\ref{el2}) implies that
$\|\,\, \cdot\,\,\|_B \leq c_2^{-1} k^{\frac{1}{2}} \|\,\, \cdot\,\,\|_2$
on $E$).\par

To prove (\ref{fin}), assume for simplicity that $A= \{1,2,\ldots,k\}$. 
For fixed
$\{g_1,g_2,\ldots,g_k\}$, and hence fixed $E$, 
$\tilde g_j = P_E g_j$, $j=k+1,k+2,\ldots,N$ are independent 
$E$--valued Gaussian
vectors with $N(0,\frac{1}{d}{\rm Id}_E)$
distribution. In particular,
${\bf E} \|\tilde g_j\|_2^2 = \frac{k}{d}$ and, by Lemma {\ref A} (ii),
$$
{\bf P}\left(\|\tilde g_j \| \geq k^{-\frac{1}{2}} \right) 
\leq \exp {\left(- \frac{k}{8} \left[\frac{k^{-1/2}}{(k/d)^{1/2}} 
\right]^2 \right)} = \exp{\left(- \frac{d}{8k} \right)}, $$
(we used the fact that $k \leq cd^{\frac{1}{2}}$, and hence 
$ k^{-\frac{1}{2}} \geq
c^{-1} (\frac{k}{d})^{\frac{1}{2}}
\geq 2 (\frac{k}{d})^{\frac{1}{2}}$).
Note a slight abuse of notation;
in fact the expectation and the probability above are conditional on 
$\{g_j \,| \,\, j \leq k\}$ To deduce (\ref{fin}), we need to know that 
$N \exp {\left(- \frac{d}{8k} \right)}$ is small. This happens e.g. when 
$k \leq \frac{d}{16\log N}$ and can be forced by the proper choice of 
$c$.\\ \pend

\noindent
{\bf Proof of Proposition \ref{square}} \,\, Clearly, it is enough to prove 
the Proposition for quotients of $l^N_1$. The variant for subspaces of 
$l^N_{\infty}$ will follow by duality.\par 
As follows from Lemma {\ref D} (i), the unit ball $B$ of a generic 
$d$--dimensional
quotient of $l^N_1$
contains a Euclidean ball $D_r$ with radius 
$r=c'\sqrt{\frac{\log {(N/d)}}{d}}$.
On the other hand, by Lemma {\ref D} (ii), $$
\left(\frac{|B|}{|D_r|}\right)^{\frac{1}{d}} \leq C'',
$$
and hence the so--called volume ratio of $B$ with respect to 
$D_r$ remainds bounded by a universal numerical constant. In particular,
(cf. [Sz-T]), this implies that, say, for $k \leq d/2$, and for a generic
$k$--dimensional subspace $G$
and some universal constant $C$,
$$
G \cap D_r \subset G \cap B \subset C(G \cap D_r), $$
which means that $G$ considered as a subspace of $B$ is $C$--Euclidean.\par

To conclude the proof we will show that
if $P$ is the orthogonal projection onto the generic $k$--dimensional 
subspace $G$ of $\Rn{N}$ and $k \geq \log N$, then
\begin{eqnarray}
\label{rzut}
P(B) \subset C_1 \sqrt{\frac{k}{d}}(G \cap D), \end{eqnarray}
where $D$ denotes the Euclidean unit ball in $\Rn{n}$. This will suffice, 
since if $k \simeq \log{N}$ and if $C$ is large enough then
$C_1\sqrt{\frac{k}{d}} \leq C c'
\sqrt{\log{\frac{N/d}{d}}}$
(remember that $N \geq d^2$ and therefore $\log N \leq 2 \log{\frac{N}{d}}$). 
Thus, (\ref{rzut}) implies $P(B) \subset CD_r$. To prove (\ref{rzut}) 
observe that for a fixed $G$ (and hence fixed $P$) and for
$j \in \{1,2,\ldots,N\}$, $Pg_j$ is a Gaussian random vector with 
distribution
$N(0, \sqrt{\frac{k}{d}}{\rm Id}_G)$.
Therefore, by Lemma {\ref A} (ii),
$$
{\bf P}\left(\|\, Pg_j \|_2 \leq
\lambda \sqrt{\frac{k}{d}} \right)
\geq 1 - \exp {(-\frac{1}{8} \lambda^2k)} $$
for $\lambda \geq 2$. Choosing
$\lambda$ sufficiently large and
using the fact that $ k \geq \log{N}$ we can obtain a similar estimate for
${\bf P} \left( \|\, Pg_j \|_2 \leq
\lambda \sqrt{\frac{k}{d}}\,|\,
j=1,2,\ldots,N \right)$. Finally, by the rotational invariance of the 
joint distribution of $g_j$'s and the Fubbini theorem we deduce that 
a ``generic'' $B$, a ``generic'' $P$ and $G$ satisfy (\ref{rzut}).
\pend

\end{document}